\documentclass[10pt,reqno]{article}

\usepackage[usenames]{color}
\usepackage{amssymb}
\usepackage{graphicx}
\usepackage{amscd}

\usepackage{graphics,amsmath,amssymb,relsize,mathalfa}
\usepackage{bm}
\usepackage{amsthm}
\usepackage{amsfonts}
\usepackage{latexsym}
\usepackage{epsf}

\newtheorem{theorem}{Theorem}[section]
\newtheorem{corollary}[theorem]{Corollary}
\newtheorem{lemma}[theorem]{Lemma}

\newtheorem{defin}[theorem]{Definition}

\newtheorem{examp}[theorem]{Example}

\newtheorem{rema}[theorem]{Remark}
\newtheorem{prob}[theorem]{Problem}

\numberwithin{equation}{section}

\newcommand{\bt}{\begin{thm}}
\newcommand{\et}{\end{thm}}
\newcommand{\bp}{\begin{proof}}
\newcommand{\ep}{\end{proof}}
\newcommand{\bprop}{\begin{prop}}
\newcommand{\eprop}{\end{prop}}
\newcommand{\bl}{\begin{lemma}}
\newcommand{\el}{\end{lemma}}
\newcommand{\bc}{\begin{corollary}}
\newcommand{\ec}{\end{corollary}}

\newcommand{\be}{\begin{enumerate}}
\newcommand{\ee}{\end{enumerate}}

\title{On a restricted linear congruence}
\author{Khodakhast Bibak \thanks{Department of Computer Science, University of Victoria, Victoria, BC, Canada V8W 3P6. Email: {\tt
\{kbibak,bmkapron,srinivas\}@uvic.ca}} \and Bruce M. Kapron \footnotemark[1] \and Venkatesh Srinivasan \footnotemark[1] \thanks{Centre for Quantum Technologies, National University of Singapore, Singapore 117543.}}

\begin{document}

\maketitle

\begin{abstract}
Let $b,n\in \mathbb{Z}$, $n\geq 1$, and ${\cal D}_1, \ldots, {\cal D}_{\tau(n)}$ be all positive divisors of $n$. For $1\leq l \leq \tau(n)$, define ${\cal C}_l:=\lbrace 1 \leqslant x\leqslant n \; : \; (x,n)={\cal D}_l\rbrace$. In this paper, by combining ideas from the finite Fourier transform of arithmetic functions and Ramanujan sums, we give a short proof for the following result: the number of solutions of the linear congruence $x_1+\cdots +x_k\equiv b \pmod{n}$, with $\kappa_{l}=|\lbrace x_1, \ldots, x_k \rbrace \cap {\cal C}_l|$, $1\leq l \leq \tau(n)$, is 
\begin{align*}
\frac{1}{n}\mathlarger{\sum}_{d\, \mid \, n}c_{d}(b)\mathlarger{\prod}_{l=1}^{\tau(n)}\left(c_{\frac{n}{{\cal D}_l}}(d)\right)^{\kappa_{l}},
\end{align*}
where $c_{d}(b)$ is a Ramanujan sum. Some special cases and other forms of this problem have been already studied by several authors. The problem has recently found very interesting applications in number theory, combinatorics, computer science, and cryptography. The above explicit formula generalizes the main results of several papers, for example, the main result of the paper by Sander and Sander [J. Number Theory {\bf 133} (2013), 705--718], one of the main results of the paper by Sander [J. Number Theory {\bf 129} (2009), 2260--2266], and also gives an equivalent formula for the main result of the paper by Sun and Yang [Int. J. Number Theory {\bf 10} (2014), 1355--1363].
\end{abstract}

\vskip .3cm
{\bf Keywords:} Restricted linear congruence; Ramanujan sum; finite Fourier transform
\vskip .3cm
{\bf 2010 Mathematics Subject Classification:} 11D79, 11P83, 11L03, 11A25, 42A16

\section{Introduction}

Recently, Bibak et al. \cite{BKSTT} gave an explicit formula for the number of solutions of the linear congruence $a_1x_1+\cdots +a_kx_k\equiv b \pmod{n}$, with $(x_i,n)=t_i$ ($1\leq i\leq k$), where $a_1,t_1,\ldots,a_k,t_k, b,n$ ($n\geq 1$) are arbitrary integers. They called these kinds of congruences {\it restricted linear congruences}. Some special cases and other forms of this problem have been already studied by several authors. The problem has recently found very interesting applications in number theory, combinatorics, computer science, and cryptography. For example, the special case of $b=0$, $a_i=1$, $t_i=\frac{n}{m_i}$, $m_i\mid n$ ($1\leq i\leq k$) is related to the {\it orbicyclic} (multivariate arithmetic) function \cite{LIS}, which has very interesting combinatorial and topological applications, in particular, in counting non-isomorphic maps on orientable surfaces \cite{LIS}. Also, the problem has been applied in studying universal hashing \cite{BKSTT2} which has many applications in computer science. Specifically, using the explicit formula for the number of solutions of the above restricted linear congruence, we designed an almost-universal hash function family and gave some applications to authentication and secrecy codes \cite{BKSTT2}.

Let $e(x)=\exp(2\pi ix)$ be the complex exponential with period 1. For integers $m$ and $n$ ($n \geq 1$) the quantity
\begin{align}\label{def1}
c_n(m) = \mathlarger{\sum}_{\substack{j=1 \\ (j,n)=1}}^{n} e\left(\frac{jm}{n}\right)
\end{align}
is called a {\it Ramanujan sum}, which is also denoted by $c(m,n)$ in the literature. From (\ref{def1}), it is clear that $c_n(-m) = c_n(m)$. Also, it is easy to see that $c_n(m)=c_n\left((m,n)\right)$, for every $m,n$.

In this paper, we prove the following theorem:

\begin{theorem} \label{thm:k var1 equi}
Let $b,n\in \mathbb{Z}$, $n\geq 1$, and ${\cal D}_1, \ldots, {\cal D}_{\tau(n)}$ be all positive divisors of $n$. For $1\leq l \leq \tau(n)$, define ${\cal C}_l:=\lbrace 1 \leqslant x\leqslant n \; : \; (x,n)={\cal D}_l\rbrace$. The number of solutions of the linear congruence $x_1+\cdots +x_k\equiv b \pmod{n}$, with $\kappa_{l}=|\lbrace x_1, \ldots, x_k \rbrace \cap {\cal C}_l|$, $1\leq l \leq \tau(n)$, is 
\begin{align}\label{thm:k var1 for equi}
\frac{1}{n}\mathlarger{\sum}_{d\, \mid \, n}c_{d}(b)\mathlarger{\prod}_{l=1}^{\tau(n)}\left(c_{\frac{n}{{\cal D}_l}}(d)\right)^{\kappa_{l}}.
\end{align}
\end{theorem}

The above theorem generalizes the main results of \cite{COH0, DIX, NV, SanSan2013}, one of the main results of \cite{San2009}, and also gives an equivalent formula for the main result of \cite{SY2014}. Note that, recently, Bibak et al. \cite{BKSTT} gave a different proof for an `equivalent' form of Theorem \ref{thm:k var1 equi}. But here we combine ideas from the finite Fourier transform of arithmetic functions and Ramanujan sums to present a new and short proof for the above theorem with the hope that its idea might be applicable to other relevant problems. In fact, as problems of this kind have many applications, especially in computer science and cryptography, having generalizations and/or new proofs and/or equivalent formulas for this problem may lead to further work. We also remark that, recently, Yang and Tang \cite{YT} considered the quadratic version of this problem in the special case of $k=2$, $a_1=a_2=1$, $t_1=t_2=1$, and posed some problems for more general cases.

\section{Proof of the theorem}

Before we proceed, we need some preliminaries. Let $r$ be a positive integer. An arithmetic function $f$ is called {\it periodic} with period $r$ (or {\it periodic} modulo $r$) if $f(m + r) = f(m)$, for every $m\in \mathbb{N}$. From (\ref{def1}), it is clear that $c_n(m)$ is a periodic function of $m$ with period 
$n$.

Let $f(m)$ be an arithmetic function with period $n$. The {\it finite Fourier transform} of $f$ is defined to be the function
\begin{align}\label{FFT1}
\widehat{f}(b)=\frac{1}{n}\mathlarger{\sum}_{m=1}^{n}f(m)e\left(\frac{-bm}{n}\right).
\end{align}
Then a Fourier representation of $f$ is obtained as (see, e.g., \cite[p. 109]{MOVA})
\begin{align}\label{FFT2}
f(m)=\mathlarger{\sum}_{b=1}^{n}\widehat{f}(b)e\left(\frac{bm}{n}\right).
\end{align}

Now we are ready to give a short proof for Theorem~\ref{thm:k var1 equi}.
\vspace{2mm}\\
\noindent\textbf{Proof of Theorem~\ref{thm:k var1 equi}.}
Suppose that $\widehat{f}_n(k,b)$ denotes the number of solutions of the linear congruence $x_1+\cdots +x_k\equiv b \pmod{n}$, with $\kappa_{l}=|\lbrace x_1, \ldots, x_k \rbrace \cap {\cal C}_l|$, $1\leq l \leq \tau(n)$. One can observe that, for every $m\in \mathbb{N}$, we have
\begin{align}\label{altprroof}
\mathlarger{\sum}_{b=1}^{n}\widehat{f}_n(k,b)e\left(\frac{bm}{n}\right) = \mathlarger{\prod}_{l=1}^{\tau(n)}\left(\mathlarger{\sum}_{x \in {\cal C}_l} e\left(\frac{mx}{n}\right)\right)^{\kappa_{l}}.
\end{align}
First, we give a short combinatorial argument to justify (\ref{altprroof}). Here the key idea is that 
$\widehat{f}_n(k,b)$ can be interpreted as the number of possible ways of writing $b$ as a sum modulo $n$  of $\kappa_{1}$ elements of ${\cal C}_1$, $\kappa_{2}$ elements of ${\cal C}_2$, $\ldots$ , 
$\kappa_{\tau(n)}$ elements of ${\cal C}_{\tau(n)}$. Now, expand the right-hand side of (\ref{altprroof}). Note that each term of this expansion has $e(\frac{m}{n})$ as a factor (compare this to the left-hand side of (\ref{altprroof})). Also note that the exponent of each term of this expansion (ignoring $m$) is just a sum of some elements of ${\cal C}_1, \ldots, {\cal C}_{\tau(n)}$, which equals $b$ ($1\leq b \leq n$). In fact, recalling the above interpretation of $\widehat{f}_n(k,b)$, we can see that in this expansion there are exactly $\widehat{f}_n(k,1)$ terms of the form $e(\frac{m}{n})$, $\widehat{f}_n(k,2)$ terms of the form $e(\frac{2m}{n})$, $\ldots$ , $\widehat{f}_n(k,n)$ terms of the form $e(m)$; that is, there are exactly $\widehat{f}_n(k,b)$ terms of the form $e(\frac{bm}{n})$, for $1\leq b \leq n$. Therefore, we get the left-hand side of (\ref{altprroof}).
\vspace*{2mm}

Putting $x_{l}'=\frac{x}{{\cal D}_l}$, $1\leq l\leq \tau(n)$, we get
\begin{align*}
\mathlarger{\sum}_{x \in {\cal C}_l} e\left(\frac{mx}{n}\right) = \mathlarger{\sum}_{\substack{x=1 \\ (x,n)={\cal D}_l}}^{n} e\left(\frac{mx}{n}\right) = \mathlarger{\sum}_{\substack{x_{l}'=1 \\ (x_{l}',n/{\cal D}_l)=1}}^{n/{\cal D}_l} e\left(\frac{mx_{l}'}{n/{\cal D}_l}\right) = c_{\frac{n}{{\cal D}_l}}(m).
\end{align*}
Therefore, 
\begin{align*}
\mathlarger{\sum}_{b=1}^{n}\widehat{f}_n(k,b)e\left(\frac{bm}{n}\right) = \mathlarger{\prod}_{l=1}^{\tau(n)}\left(c_{\frac{n}{{\cal D}_l}}(m)\right)^{\kappa_{l}}.
\end{align*}
Now, by (\ref{FFT1}) and (\ref{FFT2}), and since $c_{\frac{n}{{\cal D}_l}}\left(m\right)=c_{\frac{n}{{\cal D}_l}}\left((m,n)\right)$, we have
\begin{align*}
\widehat{f}_n(k,b) &= \frac{1}{n}\mathlarger{\sum}_{m=1}^{n}e\left(\frac{-bm}{n}\right)\mathlarger{\prod}_{l=1}^{\tau(n)}\left(c_{\frac{n}{{\cal D}_l}}(m)\right)^{\kappa_{l}}\\
&= \frac{1}{n}\mathlarger{\sum}_{d\, \mid \, n}\mathlarger{\sum}_{\substack{m=1 \\ (m, n)=d}}^{n}e\left(\frac{-bm}{n}\right)\mathlarger{\prod}_{l=1}^{\tau(n)}\left(c_{\frac{n}{{\cal D}_l}}(m)\right)^{\kappa_{l}}\\
&= \frac{1}{n}\mathlarger{\sum}_{d\, \mid \, n}\mathlarger{\sum}_{\substack{m=1 \\ (m, n)=d}}^{n}e\left(\frac{-bm}{n}\right)\mathlarger{\prod}_{l=1}^{\tau(n)}\left(c_{\frac{n}{{\cal D}_l}}(d)\right)^{\kappa_{l}}\\
&{\stackrel{m'=m/d}{=}} \frac{1}{n}\mathlarger{\sum}_{d\, \mid \, n}\mathlarger{\sum}_{\substack{m'=1 \\ (m', n/d)=1}}^{n/d}e\left(\frac{-bm'}{n/d}\right)\mathlarger{\prod}_{l=1}^{\tau(n)}\left(c_{\frac{n}{{\cal D}_l}}(d)\right)^{\kappa_{l}}\\
&= \frac{1}{n}\mathlarger{\sum}_{d\, \mid \, n}c_{n/d}(-b)\mathlarger{\prod}_{l=1}^{\tau(n)}\left(c_{\frac{n}{{\cal D}_l}}(d)\right)^{\kappa_{l}}\\
&= \frac{1}{n}\mathlarger{\sum}_{d\, \mid \, n}c_{n/d}(b)\mathlarger{\prod}_{l=1}^{\tau(n)}\left(c_{\frac{n}{{\cal D}_l}}(d)\right)^{\kappa_{l}}
= \frac{1}{n}\mathlarger{\sum}_{d\, \mid \, n}c_{d}(b)\mathlarger{\prod}_{l=1}^{\tau(n)}\left(c_{\frac{n}{{\cal D}_l}}(d)\right)^{\kappa_{l}}.
\end{align*} \hfill $\Box$

\section*{Acknowledgements}

The authors would like to thank the anonymous referees for helpful comments. During the preparation of this work the first author was supported by a Fellowship from the University of Victoria (UVic Fellowship).

\end{document}